\documentclass[11pt,reqno]{amsart}
\usepackage{amsmath,amsthm,amssymb,verbatim,enumerate,ifthen}
\usepackage[mathscr]{eucal}
\usepackage[utf8]{inputenc}
\def\N{\mathbb{N}}
\def\R{\mathbb{R}}

\def\B{\mathscr{B}}
\def\CM{\mathscr{CM}}
\def\CP{\mathscr{CP}}

\newtheorem{theorem}{Theorem}
\newtheorem*{theorem*}{Theorem}
\def\Thm#1#2{\ifthenelse{\equal{#1}{*}}{\begin{theorem*}#2\end{theorem*}}
  {\begin{theorem}\label{T#1}#2\end{theorem}}}
\newtheorem{Atheorem}{Theorem}

\def\thm#1{Theorem~\ref{T#1}}

\newtheorem{proposition}[theorem]{Proposition}
\newtheorem*{proposition*}{Proposition}
\def\Prp#1#2{\ifthenelse{\equal{#1}{*}}{\begin{proposition*}#2\end{proposition*}}
             {\begin{proposition}\label{P#1}#2\end{proposition}}}

\newtheorem{corollary}[theorem]{Corollary}
\newtheorem*{corollary*}{Corollary}
\def\Cor#1#2{\ifthenelse{\equal{#1}{*}}{\begin{corollary*}#2\end{corollary*}}
             {\begin{corollary}\label{C#1}#2\end{corollary}}}
\def\cor#1{Corollary~\ref{C#1}}

\newtheorem{lemma}[theorem]{Lemma}
\newtheorem*{lemma*}{Lemma}
\def\Lem#1#2{\ifthenelse{\equal{#1}{*}}{\begin{lemma*}#2\end{lemma*}}
             {\begin{lemma}\label{L#1}#2\end{lemma}}}
\def\lem#1{Lemma~\ref{L#1}}

\newtheorem{example}[theorem]{Example}
\newtheorem*{example*}{Example}
\def\Exa#1#2{\ifthenelse{\equal{#1}{*}}{\begin{example*}\rm #2\end{example*}}
             {\begin{example}\label{Ex#1}\rm #2\end{example}}}

\newtheorem{problem}[theorem]{Problem}

\theoremstyle{definition}
\newtheorem{definition}[theorem]{Definition}

\newtheorem{remark}[theorem]{Remark}
\newtheorem*{remark*}{Remark}
\def\Rem#1#2{\ifthenelse{\equal{#1}{*}}{\begin{remark*}\rm #2\end{remark*}}
             {\begin{remark}\label{R#1}\rm #2\end{remark}}}

\def\eq#1{{\rm(\ref{E#1})}}
\def\Eq#1#2{\ifthenelse{\equal{#1}{*}}
  {\begin{equation*}\begin{aligned}#2\end{aligned}\end{equation*}}
  {\begin{equation}\begin{aligned}\label{E#1}#2\end{aligned}\end{equation}}}

\long\def\comment#1{}

\begin{document}
\begin{flushright}
\end{flushright}
\vspace{5mm}

\date{\today}
\title[On the equality problem of two-variable Bajraktarevi\'c means]{On the equality problem of two-variable Bajraktarevi\'c means under first-order differentiability assumptions}

\author[Zs. P\'ales]{Zsolt P\'ales}
\address[Zs. P\'ales]{Institute of Mathematics, University of Debrecen, H-4002 Debrecen, Pf.\ 400, Hungary}
\email{pales@science.unideb.hu}

\author[A. Zakaria]{Amr Zakaria}
\address[A. Zakaria]{Department of Mathematics, Faculty of Education, Ain Shams University, Cairo 11341, Egypt}
\email{amr.zakaria@edu.asu.edu.eg}

\subjclass[2010]{39B22, 39B52}
\keywords{generalized quasiarithmetic mean; equality problem; functional equation, regularity problem}


\thanks{The research of the first author was supported by the K-134191 NKFIH Grant and the 2019-2.1.11-TÉT-2019-00049, the EFOP-3.6.1-16-2016-00022 and the EFOP-3.6.2-16-2017-00015 projects. The last two projects are co-financed by the European Union and the European Social Fund. The research of the second author was supported by the Bilateral State Scholarship of the Tempus Public Foundation of Hungary BE AK 2020-2021/157507.}

\begin{abstract}
The equality problem of the two-variable Bajraktarevi\'c means can be expressed as the functional equation
$$
  \left(\frac{f}{g}\right)^{-1}\bigg(\frac{f(x)+f(y)}{g(x)+g(y)}\bigg)
   =\left(\frac{h}{k}\right)^{-1}\bigg(\frac{h(x)+h(y)}{k(x)+k(y)}\bigg)\qquad(x,y\in I),
$$
where $I$ is a nonempty open real interval, $f,g,h,k:I\to\R$ are continuous functions, $g$, $k$ are positive and $f/g$, $h/k$ are strictly monotone.
This functional equation, for the first time, was solved by Losonczi in 1999 under 6th-order continuous differentiability assumptions. Additional and new characterizations of this equality problem have been found recently by Losonczi, P\'ales and Zakaria under the same regularity assumptions in 2021. In this paper it is shown that the same conclusion can be obtained under substantially weaker regularity conditions, namely, assuming only first-order differentiability.
\end{abstract}

\maketitle

\section{\bf Introduction}

In this paper $I$ will stand for a nonempty open interval and the classes of continuous strictly monotone and continuous positive functions defined on $I$ will be denoted by $\CM(I)$ and $\CP(I)$, respectively.

Given two functions $f,g:I\to\R$ with $g\in\CP(I)$ and $f/g\in\CM(I)$, the \emph{two-variable Bajraktarevi\'c mean} $\B_{f,g}:I^2\to I$ is defined by
\Eq{*}{
  \B_{f,g}(x,y)
      :=\left(\frac{f}{g}\right)^{-1}\bigg(\frac{f(x)+f(y)}{g(x)+g(y)}\bigg)
   \qquad(x,y\in I).
}
The pair $(f,g)$ is called the \emph{generator of $\B_{f,g}$}.
This notion was introduced, in a somewhat different form and for more than two variables, by Bajraktarevi\'c \cite{Baj58}, \cite{Baj63} in the late 50s. By taking $g\equiv1$, we can see that this class of means includes quasiarithmetic means (cf.\ \cite{HarLitPol34}). 

One of the most basic questions related to this class of means is their \emph{equality problem}, which is to find the solutions $f,g,h,k:I\to\R$ of the functional equation
\Eq{1}{
   \B_{f,g}(x,y)=\B_{h,k}(x,y)\qquad(x,y\in I).
}
It is not difficult to prove that if there exist four real constants $a,b,c,d$ with $ad-bc\neq0$ such that
\Eq{hk}{
  h=af+bg \qquad\mbox{and}\qquad k=cf+dg
}
hold on $I$, then the means $\B_{f,g}$ and $\B_{h.k}$ are identical. Motivated by this observation, if the equalities \eq{hk} are satisfied by some constants $a,b,c,d$ with $ad-bc\neq0$, then we say that \emph{the pair $(f,g)$ is equivalent to $(h,k)$} and we denote this property by $(f,g)\sim(h,k)$. For the equality of two Bajraktarevi\'c means with more than two variables the equivalence of their generators is not only sufficient but it is also necessary (cf.\ \cite{Baj58}, \cite{Baj63}, \cite{DarLos70}, \cite{PalZak20b}, \cite{PalZak20a}). Surprisingly, in the setting of two-variable means, the situation is completely different. This case was investigated by Losonczi \cite{Los99}, \cite{Los06b} who, assuming 6 times continuous differentiability, discovered that beyond the equality via equivalence of the generators, there are further 32 possibilities.

This result has been revisited by Losonczi, P\'ales, and Zakaria in \cite{LosPalZak21}, where several new characterizations of the equality of two-variable Bajraktarevi\'c means \eq{1} have been obtained. To recall this result, we need to introduce some further notations. If $f,g:I\to\R$ are differentiable functions, then their Wronskian $W_{f,g}$ is defined as
\Eq{*}{
  W_{f,g}:=f'g-fg'.
}
For a real parameter $\alpha$, let us define the functions $S_\alpha,C_\alpha:\R\to\R$ as follows:
\Eq{*}{
  S_\alpha(t)
  &:=\sum_{k=0}^\infty\frac{\alpha^k t^{2k+1}}{(2k+1)!}
  =\begin{cases}
  \dfrac{\sin(\sqrt{-\alpha}t)}{\sqrt{-\alpha}}  
  &\text{ if } \alpha<0,\\ 
  t &\text{ if } \alpha=0,\\
  \dfrac{\sinh(\sqrt{\alpha}t)}{\sqrt{\alpha}}  
  &\text{ if } \alpha>0, 
  \end{cases}\\
  C_\alpha(t)
  &:=\sum_{k=0}^\infty\frac{\alpha^k t^{2k}}{(2k)!}
  =\begin{cases}
  \cos(\sqrt{-\alpha}t)  
  &\text{ if } \alpha<0,\\[2mm] 
  1 &\text{ if } \alpha=0,\\[2mm]
  \cosh(\sqrt{\alpha}t) 
  &\text{ if } \alpha>0. 
  \end{cases}
}
We note that these functions form a fundamental system of solutions of the second-order homogeneous linear differential equation $Y''=\alpha Y$.

\Thm{EBM}{Let $f,g,h,k:I\to\R$ be 6 times continuously differentiable functions such that $g,k\in\CP(I)$ and the derivatives $(f/g)'$, $(h/k)'$ do not vanish in $I$. Then the following assertions are equivalent:
\begin{enumerate}[(i)]
 \item The means $\B_{f,g}$ and $\B_{h.k}$ satisfy the functional equation \eq{1}.
 \item Either $(f,g)\sim(h,k)$ or there exist real constants $\alpha$, $\beta$ and $\gamma$ such that
 \Eq{*}{
    W_{f',g'}=\alpha\big(W_{f,g}\big)^3, \qquad
    W_{h',k'}=\beta\big(W_{h,k}\big)^3, \qquad\mbox{and}\qquad
    W_{h,k}=\gamma W_{f,g}
 }
 hold on $I$.
 \item Either $(f,g)\sim(h,k)$ or there exist real constants $a,b,c,A,B,C,\gamma$ such that
 \Eq{*}{
   af^2+bfg+cg^2=1, \qquad Ah^2+Bhk+Ck^2=1,
   \qquad\mbox{and}\qquad W_{h,k}=\gamma W_{f,g}
 }
 hold on $I$.
 \item Either $(f,g)\sim(h,k)$ or there exist two real polynomials $P$ and $Q$ of at most second degree which are positive on the ranges of $f/g$ and $h/k$, respectively, and there exist real constants $\gamma$ and $\delta$ such that
 \Eq{*}{
   g=\frac{1}{\sqrt{P}}\circ\frac{f}{g}, \qquad
   k=\frac{1}{\sqrt{Q}}\circ\frac{h}{k}, \qquad\mbox{and}\qquad
   \bigg(\int\frac{1}{Q}\bigg)\circ\frac{h}{k}
   =\gamma\bigg(\int\frac{1}{P}\bigg)\circ\frac{f}{g}+\delta.
 }
 \item Either $(f,g)\sim(h,k)$ or there exist a continuous strictly monotone function $w:I\to\R$ and real constants $\alpha$ and $\beta$ such that
 \Eq{*}{
   (f,g)\sim(S_\alpha\circ w,C_\alpha\circ w) \qquad\mbox{and}\qquad
   (h,k)\sim(S_\beta\circ w,C_\beta\circ w).
 }
 \item Either $(f,g)\sim(h,k)$ or $\B_{f,g}=\B_{w,1}=\B_{h.k}$ holds on $I^2$ with $w:=\int W_{f,g}$.
 \item Either $(f,g)\sim(h,k)$ or there exists a continuous strictly monotone function $w:I\to\R$ such that  $\B_{f,g}=\B_{w,1}=\B_{h.k}$ holds on $I^2$.
\end{enumerate}}

In fact, assertion (ii) of the above theorem is the corrected form of assertion (iv) of \cite[Theorem 15]{LosPalZak21}, where, unfortunately, the condition $(f,g)\sim(h,k)$ was missing.

The equality problem of means in various classes of two-variable means has been solved. We refer here to Losonczi's works \cite{Los99}, \cite{Los00a}, \cite{Los02a}, \cite{Los03a}, \cite{Los06b} where the equality of two-variable means is characterized in several settings. A key idea in these papers, under appropriate order differentiability assumptions, is to calculate and then to compare the partial derivatives of the means at points of the form $(x,x)$, where $x\in I$. A similar problem, the mixed equality problem of quasiarithmetic and Lagrangian means was solved by P\'ales \cite{Pal11}. The equality problem of generalized Bajraktarevi\'c means was solved by Gr\"unwald and P\'ales in \cite{GruPal20}.

The main goal of this paper is to solve one of the open problems listed in the paper \cite{Pal02a}. We show that conditions (i) and (iii)--(vii) of \thm{EBM} are equivalent to each other under first-order continuous differentiability assumptions and thus, we substantially weaken the regularity requirements of \thm{EBM}. The main difficulty to prove our result comes the fact that the validity of equality \eq{1} does not imply the differentiability properties of the unknown functions $f,g,h,k$. Indeed, \eq{1} can hold only if $(f,g)\sim(h,k)$. On the other hand, if \eq{1} holds but  $(f,g)\not\sim(h,k)$, then the functions $p,q,\varphi,\psi:J\to\R$ which are constructed in \lem{fghk} are infinitely many times differentiable. 

The pairwise equivalence of conditions (i) and (iii)--(vii) suggests that statement of \thm{EBM} (without assertion (ii) could be true even if we do not assume first-order continuous differentiability.

\section{Reduction of the equality problem and auxiliary results}

In what follows, we reduce the equality problem of two-variable Bajraktarevi\'c means to the equality of a two-variable Bajraktarevi\'c mean and an arithmetic mean weighted by a weight function.

\Lem{fghk}{Let $f,g,h,k:I\to\R$ be continuous functions such that $g,k\in\CP(I)$, $f/g,h/k\in\CM(I)$. Then the equality \eq{1} holds if and only if 
\Eq{2}{
  \B_{q,p}(u,v)=\B_{\varphi,\psi}(u,v)\qquad(u,v\in J),
}
where
\Eq{*}{
  p:=k\circ\bigg(\frac{h}{k}\bigg)^{-1}, \qquad 
  q(u):=p(u)u, \qquad \varphi:=f\circ\bigg(\frac{h}{k}\bigg)^{-1}, 
  \qquad\psi:=g\circ\bigg(\frac{h}{k}\bigg)^{-1},
  \qquad J:=\big(\tfrac{h}{k}\big)(I).
}}

\begin{proof}
Equation \eq{1} is equivalent to the functional equation
\Eq{fghk}{
  \bigg(\frac{f}{g}\bigg)^{-1}\bigg(\frac{f(x)+f(y)}{g(x)+g(y)}\bigg)
  =\bigg(\frac{h}{k}\bigg)^{-1}\bigg(\frac{h(x)+h(y)}{k(x)+k(y)}\bigg)
  \qquad(x,y\in I).
}
Then, with the substitutions $x:=\big(\frac{h}{k}\big)^{-1}(u)$ and $y:=\big(\frac{h}{k}\big)^{-1}(v)$, \eq{fghk} implies
\Eq{rep}{
  \bigg(\frac{\varphi}{\psi}\bigg)^{-1}
  \bigg(\frac{\varphi(u)+\varphi(v)}{\psi(u)+\psi(v)}\bigg)
  =\bigg(\frac{h}{k}\bigg)\circ\bigg(\frac{f}{g}\bigg)^{-1}\bigg(\frac{f(x)+f(y)}{g(x)+g(y)}\bigg)
  =\frac{h(x)+h(y)}{k(x)+k(y)}=\frac{p(u)u+p(v)v}{p(u)+p(v)}
}
for all $u,v\in J:=\big(\tfrac{h}{k}\big)(I)$. This proves that \eq{2} is valid.

The proof of the reversed implication is similar.
\end{proof}

\Lem{1}{Let $n\in\N$ and $\varphi,\psi:J\to\R$ be $n$-times continuously differentiable functions such that $\psi\in\CP(J)$ and $\big(\frac{\varphi}{\psi}\big)'$ is nowhere vanishing. Then $\B_{\varphi,\psi}$ is $n$-times continuously differentiable on $J^2$.}

\begin{proof}
Due to the assumptions, the function $\big(\frac{\varphi}{\psi}\big)^{-1}$ is $n$-times continuously differentiable on $\frac{\varphi}{\psi}(J)$. By the chain rule, this implies that the mean $\B_{\varphi,\psi}$ is $n$-times continuously differentiable on $J^2$.
\end{proof}

\Lem{3}{Let $n\in\N$ and $\varphi,\psi:J\to\R$ be $n$-times continuously differentiable functions such that $\psi\in\CP(J)$ and $\big(\frac{\varphi}{\psi}\big)'$ is nowhere vanishing. Let $p\in\CP(J)$ and define $q:J\to\R$ by $q(u):=p(u)u$. If \eq{2} holds, then $p$ is $n$-times continuously differentiable.}

\begin{proof} 
To prove the $n$-times continuous differentiability of $p$ at an arbitrarily fixed element $u_0\in J$, let $v_0\in J$ be another fixed element which is distinct from $u_0$. Then, due to the strict mean property of Bajraktarevi\'c means, it follows that $\B_{\varphi,\psi}(u_0,v_0)\neq u_0$. Hence, for some neighborhood $U_0$ of $u_0$, we have that $\B_{\varphi,\psi}(u,v_0)\neq u$ hold for all $u\in U_0$. By the equality \eq{2}, for $u\in U_0$, we get
\Eq{v0}{
	p(u)
	=p(v_0)\cdot\frac{v_0-\B_{\varphi,\psi}(u,v_0)}{\B_{\varphi,\psi}(u,v_0)-u}.
}
Using this formula, and applying \lem{1}, we can see that $p$ is $n$-times continuously differentiable on $U_0$, in particular, at $u_0$. 
\end{proof}

\section{Main Results}

First, we are going to establish various consequences of the equality \eq{2}.

\Lem{4}{Let $p:J\to\R$ and let $\varphi,\psi:J\to\R$ be continuously differentiable functions such that $\psi,p\in\CP(J)$ and $\big(\frac{\varphi}{\psi}\big)'$ is nowhere vanishing. Define $q:J\to\R$ by $q(u):=p(u)u$. If \eq{2} holds, then $p$ is continuously differentiable and, for all $(u,v)\in J^2$, 
\Eq{prod}{
	&\Big(\varphi'(u)(\psi(u)+\psi(v))-(\varphi(u)+\varphi(v))\psi'(u)\Big)
	\Big(p'(v)p(u)(v-u)+p(v)(p(u)+p(v))\Big)\\
	&=\Big(\varphi'(v)(\psi(u)+\psi(v))-(\varphi(u)+\varphi(v))\psi'(v)\Big)
	\Big(p'(u)p(v)(u-v)+p(u)(p(u)+p(v))\Big).
}
Furthermore, there exists a nonzero constant $\gamma\in\R$ such that
\Eq{cp}{
	W_{\varphi,\psi}=\gamma p^2.
}}

\begin{proof}
By \lem{1} and \lem{3}, $p$ and $\B_{\varphi,\psi}$ are continuously differentiable on $J$ and $J^2$, respectively. Differentiating the equation \eq{2} with respect to the first and second variables, for all $(u,v)\in J^2$, we obtain
\Eq{*}{
  \bigg(\bigg(\frac{\varphi}{\psi}\bigg)^{-1}\bigg)'&
  \bigg(\frac{\varphi(u)+\varphi(v)}{\psi(u)+\psi(v)}\bigg)
  \frac{\varphi'(u)(\psi(u)+\psi(v))-(\varphi(u)+\varphi(v))\psi'(u)}{(\psi(u)+\psi(v))^2}\\
  &=\frac{p'(u)p(v)(u-v)+p(u)(p(u)+p(v))}{(p(u)+p(v))^2},\\
  \bigg(\bigg(\frac{\varphi}{\psi}\bigg)^{-1}\bigg)'&
  \bigg(\frac{\varphi(u)+\varphi(v)}{\psi(u)+\psi(v)}\bigg)
  \frac{\varphi'(v)(\psi(u)+\psi(v))-(\varphi(u)+\varphi(v))\psi'(v)}{(\psi(u)+\psi(v))^2}\\
  &=\frac{p'(v)p(u)(v-u)+p(v)(p(u)+p(v))}{(p(u)+p(v))^2}.
}
These equalities easily yield that \eq{prod} is true for all $(u,v)\in J^2$. 

Now, using the definition of $W_{\varphi,\psi}$, the equality \eq{prod} can be rewritten as
\Eq{*}{
	&\Big(2W_{\varphi,\psi}(u)+\varphi'(u)(\psi(v)-\psi(u))-(\varphi(v)-\varphi(u))\psi'(u)\Big)
	\Big(p'(v)p(u)(v-u)+p(v)(p(u)+p(v))\Big)\\
	&=\Big(2W_{\varphi,\psi}(v)+\varphi'(v)(\psi(u)-\psi(v))-(\varphi(u)-\varphi(v))\psi'(v)\Big)
	\Big(p'(u)p(v)(u-v)+p(u)(p(u)+p(v))\Big).
}
After rearranging this equality, we arrive at
\Eq{*}{
  \frac{\big(W_{\varphi,\psi}/p\big)(v)-\big(W_{\varphi,\psi}/p\big)(u)}{v-u}
  &=W_{\varphi,\psi}(u)\frac{p'(v)}{p(v)(p(u)+p(v))}
  	+W_{\varphi,\psi}(v)\frac{p'(u)}{p(u)(p(u)+p(v))}\\
  &\quad-\frac{\varphi'(v)(\psi(u)-\psi(v))-(\varphi(u)-\varphi(v))\psi'(v)}{2(v-u)}
  \Big(\frac{p'(u)(u-v)}{p(u)(p(u)+p(v))}+\frac{1}{p(v)}\Big)\\
  &\quad+\frac{\varphi'(u)(\psi(v)-\psi(u))-(\varphi(v)-\varphi(u))\psi'(u)}{2(v-u)}
  \Big(\frac{p'(v)(v-u)}{p(v)(p(u)+p(v))}+\frac{1}{p(u)}\Big).
}
Observe now that each term on the right hand side of this equality has a limit as $v$ tends to $u$. Therefore, the limit of the left hand side must exist, showing that the function $w:=W_{\varphi,\psi}/p$ is differentiable at $J$. Upon taking the limit as $v\to u$, we obtain
\Eq{*}{
  w'(u)
  =\frac{p'}{p}(u)w(u),
}
which implies that $w/p$ is constant on $J$. Consequently, there exists $\gamma\in\R$ such that \eq{cp} holds on $J$. The constant $\gamma$ cannot be zero, otherwise $W_{\varphi,\psi}$ is identically zero, which contradicts that the derivative $\big(\frac{\varphi}{\psi}\big)'$ is nowhere vanishing. 
\end{proof}

\Lem{Reg}{Let $\varphi,\psi,p:J\to\R$ be continuously differentiable functions such that $\psi,p\in\CP(J)$ and $\big(\frac{\varphi}{\psi}\big)'$ is nowhere vanishing. Assume that the equality \eq{prod} is satisfied for all $(u,v)\in J^2$. Then there exists an open (possibly empty) subset $H\subseteq J$ such that $\varphi,\psi,p$ are infinitely many times differentiable on $H$ and 
\Eq{Wid}{
\gamma W_{q',p'}-W_{\varphi',\psi'}=\gamma\frac{r''}rp^2
}
holds on $H$, where $r:=\psi/p$. If $K:=J\setminus H$ is not a singleton then $r$ is twice continuously differentiable on $J$ and $r''$ vanishes on $K$. Finally, if $K=\{u_0\}$ is valid for some $u_0\in J$, then $r$ is twice continuously differentiable on $J\setminus\{u_0\}$. If, additionally, $r''$ has finite left and right limits at $u_0$, then $r$ is also twice continuously differentiable at $u_0$ and $r''(u_0)=0$.}

\begin{proof} Assume that \eq{prod} holds for all $(u,v)\in J^2$. In view of \lem{4}, it follows that there exists a nonzero constant $\gamma\in\R$ such that \eq{cp} holds in $J$. Therefore,
\Eq{Wd}{
  \Big(\frac{\varphi}{\psi}\Big)'
  =\frac{W_{\varphi,\psi}}{\psi^2}
  =\frac{\gamma p^2}{\psi^2}.
}
Thus, for all $(u,v)\in J^2$, we get 
\Eq{*}{
  \varphi'(u)(\psi(u)+\psi(v))-(\varphi(u)+\varphi(v))\psi'(u)
  &=\Big(\frac{\varphi}{\psi}\Big)'(u)\psi(u)(\psi(u)+\psi(v))+\psi'(u)\psi(v)\int_v^u\Big(\frac{\varphi}{\psi}\Big)'\\
  &=\gamma\bigg(\frac{p^2(u)(\psi(u)+\psi(v))}{\psi(u)}+\psi'(u)\psi(v)\int_v^u\frac{p^2}{\psi^2}\bigg).
}
Using this formula, the equality \eq{prod} can equivalently be rewritten as
\Eq{*}{
	&\bigg(\frac{p^2(u)(\psi(u)+\psi(v))}{\psi(u)}+\psi'(u)\psi(v)\int_v^u \frac{p^2}{\psi^2}\bigg)
	\bigg(\frac{p'(v)(v-u)}{p(v)}+\frac{p(u)+p(v)}{p(u)}\bigg)\\
	&\qquad=\bigg(\frac{p^2(v)(\psi(u)+\psi(v))}{\psi(v)}+\psi'(v)\psi(u)\int_u^v \frac{p^2}{\psi^2}\bigg)
	\bigg(\frac{p'(u)(u-v)}{p(u)}+\frac{p(u)+p(v)}{p(v)}\bigg).
}
Therefore, for all $(u,v)\in J^2$ with $u\neq v$,
this is an equation of the following form:
\Eq{ABC}{
  A(u,v)\psi'(u)+B(u,v)p'(u)=C(u,v),
}
where
\Eq{*}{
  A(u,v)&:=\bigg(\frac{(p(u)+p(v))\psi(v)}{p(u)(u-v)}-\frac{p'(v)\psi(v)}{p(v)}\bigg)\int_u^v \frac{p^2}{\psi^2},\\
  B(u,v)&:=\frac{p^2(v)(\psi(u)+\psi(v))}{p(u)\psi(v)}+\frac{\psi(u)\psi'(v)}{p(u)}\int_u^v \frac{p^2}{\psi^2},\\
  C(u,v)&:=\bigg(\frac{p(u)}{\psi(u)}-\frac{p(v)}{\psi(v)}\bigg)\frac{(\psi(u)+\psi(v))(p(u)+p(v))}{u-v}\\
  &\qquad-\frac{p'(v)p^2(u)(\psi(u)+\psi(v))}{p(v)\psi(u)}
  -\frac{(p(u)+p(v))\psi(u)\psi'(v)}{p(v)(u-v)}\int_u^v \frac{p^2}{\psi^2}.
}
For a fixed $u\in J$, consider now the curve $\Gamma(u)$ defined by
\Eq{*}{
  \Gamma(u):=\big\{(A(u,v),B(u,v))\colon v\in J\setminus\{u\}\big\}\subseteq\R^2
}
and define
\Eq{*}{
  K:=\{u\in J\colon \dim\Gamma(u)=1\}\qquad\mbox{and}\qquad
  H:=\{u\in J\colon \dim\Gamma(u)=2\}.
}
Obviously, $\{H,K\}$ form a partition of $J$, i.e., they are disjoint and their union equals $J$. First, we are going to show that $H$ is open and $p,\psi,\varphi$ are infinitely many times differentiable over $H$. 

Let $u_0\in H$ be fixed. Then there exist two independent elements of $\Gamma(u_0)$, i.e., there exist $v_1,v_2\in J\setminus\{u_0\}$ such that 
\Eq{*}{
A(u_0,v_1)B(u_0,v_2)\neq A(u_0,v_2)B(u_0,v_1).
}
By the continuity of $A$ and $B$, this implies that there exists a neighborhood $U_0\subseteq J\setminus\{v_1,v_2\}$ of $u_0$ such that
\Eq{*}{
A(u,v_1)B(u,v_2)\neq A(u,v_2)B(u,v_1),
}
which implies that $\dim\Gamma(u)=2$ for all $u\in U_0$.  Therefore, $U_0\subseteq H$ proving that $H$ is open. Substituting $v=v_i$ into \eq{ABC} and solving the system of equations so obtained with respect to the pair $(\psi'(u),p'(u))$, we get
\Eq{*}{
  \psi'(u)=\frac{C(u,v_1)B(u,v_2)-C(u,v_2)B(u,v_1)}{A(u,v_1)B(u,v_2)-A(u,v_2)B(u,v_1)},\qquad
  p'(u)=\frac{A(u,v_1)C(u,v_2)-A(u,v_2)C(u,v_1)}{A(u,v_1)B(u,v_2)-A(u,v_2)B(u,v_1)}
}
for all $u\in U_0$. Assume that we have proved that $\psi$ and $p$ are $n$-times continuously differentiable on $U_0$. Then $A(\cdot,v_i)$, $B(\cdot,v_i)$, and $C(\cdot,v_i)$ are also $n$-times continuously differentiable on $U_0$. Thus, the above equalities imply that $\psi'$ and $p'$ are $n$-times continuously differentiable on $U_0$, too. That is, $\psi$ and $p$ are $(n+1)$-times continuously differentiable on $U_0$ and hence on $H$. Having proved the infinitely many times differentiability of $\psi$ and $p$ on $H$, the equality \eq{Wd} implies that $\varphi$ is also infinitely many times differentiable on $H$.

We have that $q(u)=up(u)$. Therefore, for all $u\in H$,
\Eq{*}{
W_{q,p}(u)
&=q'(u)p(u)-q(u)p'(u)
=(p(u)+up'(u))p(u)-up(u)p'(u)=p^2(u),\\
W_{q',p'}(u)
&=q''(u)p'(u)-q'(u)p''(u)
=(2p'(u)+up''(u))p'(u)-(p(u)+up'(u))p''(u)
=(2p'^2-pp'')(u).\\
}
On the other hand,
\Eq{*}{
\gamma p^2=W_{\varphi,\psi}=\varphi'\psi-\varphi\psi'
\qquad\mbox{and}\qquad 
2\gamma p'p =\big(W_{\varphi,\psi}\big)'=\varphi''\psi-\varphi\psi''.
}
Hence, 
\Eq{*}{
\varphi'=\frac{\gamma p^2+\varphi\psi'}{\psi}\qquad\mbox{and}\qquad\varphi''=\frac{2\gamma p'p+\varphi\psi''}{\psi}.
}
Therefore, using also that $\psi=rp$, we get
\Eq{*}{
W_{\varphi',\psi'}
&=\varphi''\psi'-\varphi'\psi''
=\frac{2\gamma p'p+\varphi\psi''}{\psi}\psi'-\frac{\gamma p^2+\varphi\psi'}{\psi}\psi''
=\frac{2\gamma p'p}{\psi}\psi'-\frac{\gamma p^2}{\psi}\psi''\\
&=\frac{\gamma}r(2p'\psi'-p\psi'')
=\frac{\gamma}r\big(2p'(r'p+rp')-p(r''p+2r'p'+rp'')\big)\\
&=\gamma\Big(2p'^2-p''p-\frac{r''}rp^2\Big)
=\gamma\Big(W_{q',p'}-\frac{r''}rp^2\Big).
}
This implies that the identity \eq{Wid} holds on $H$.

Let now $u\in K$ be arbitrarily fixed. Then there exists $(\lambda(u),\mu(u))\neq(0,0)$ such that
\Eq{*}{
  \lambda(u)A(u,v)+\mu(u) B(u,v)=0 \qquad(v\in J\setminus\{u\}).
}
This gives, for all $v\in J\setminus\{u\}$ that
\Eq{ab}{
  &\lambda(u)\bigg(\frac{(p(u)+p(v))\psi(v)}{p(u)(u-v)}-\frac{p'(v)\psi(v)}{p(v)}\bigg)\int_u^v \frac{p^2}{\psi^2}\\
  &+\mu(u)\bigg(\frac{p^2(v)(\psi(u)+\psi(v))}{p(u)\psi(v)}+\frac{\psi(u)\psi'(v)}{p(u)}\int_u^v \frac{p^2}{\psi^2}\bigg)=0.
}
Taking the limit $v\to u$, we get
\Eq{*}{
  -2\lambda(u)\frac{p^2(u)}{\psi(u)}
  +2\mu(u)p(u)=0,
}
which implies that $\lambda(u)p(u)=\mu(u)\psi(u)$. Therefore, $\lambda(u)\neq0\neq\mu(u)$ and thus \eq{ab}, for all $v\in J\setminus\{u\}$, can be rewritten as
\Eq{*}{
  \bigg(\frac{p(u)+p(v)}{p(u)(v-u)}+\frac{p'(v)}{p(v)}-\frac{\psi'(v)}{\psi(v)}\bigg)\int_u^v \frac{p^2}{\psi^2}
  =p^2(v)\frac{\psi(u)+\psi(v)}{\psi(u)\psi(v)^2}.
}
Let us substitute $\psi=rp$ into this equality. Then we obtain
\Eq{rp}{
  \bigg(\frac{p(u)+p(v)}{p(u)(v-u)}-\frac{r'(v)}{r(v)}\bigg)\int_u^v \frac{1}{r^2}
  =\frac{p(u)r(u)+p(v)r(v)}{p(u)r(u)r(v)^2}.
}
Solving this equation with respect to $r'(v)$ (when $u\neq v\in J$), we get
\Eq{*}{
  r'(v)=\frac{r(v)(p(u)+p(v))}{p(u)(v-u)}
  -\bigg(\frac{1}{r(v)}+\frac{p(v)}{p(u)r(u)}\bigg)\bigg(\int_u^v \frac{1}{r^2}\bigg)^{-1},
}
from where we can see that $r'$ is continuously differentiable at every $v\in J$ which is different from $u$. The equality \eq{rp} implies that 
\Eq{*}{
  \int_u^v \frac{1}{r^2}
  =\frac{(v-u)(p(u)r(u)+p(v)r(v))}
  {r(v)r(u)\big(r(v)(p(u)+p(v))-p(u)r'(v)(v-u)\big)} \qquad (v\in J\setminus\{u\}).
}
Differentiating this equality with respect to $v\in J\setminus\{u\}$, we get
\Eq{*}{
  &\big(r(v)r(u)(p(u)+p(v))-p(u)r'(v)r(u)(v-u)\big)^2\\
  &=\big(p(u)r(u)+p(v)r(v)+(v-u)(p'(v)r(v)+p(v)r'(v))\big)\\&\hspace{1.5cm}\cdot r(v)r(u)\big(r(v)(p(u)+p(v))-p(u)r'(v)(v-u)\big)\\
  &\quad-(v-u)(p(u)r(u)+p(v)r(v))\\
  &\hspace{1.5cm}\cdot\big(r'(v)r(v)r(u)(p(u)+2p(v))+r(v)^2r(u)p'(v)-p(u)(r''(v)r(v)+r'(v)^2)r(u)(v-u)\big),
}
whence we get
\Eq{r''}{
  r''(v)=\frac{r(v)\big(p(v)(p(u)+p(v))
  +(v-u)p'(v)p(u)\big)}{(p(u)r(u)+p(v)r(v))p(u)}
  \cdot\frac{r(u)-r(v)+r'(v)(v-u)}{(v-u)^2}.
}
By the Cauchy Mean Value Theorem, there exists a point $w=w(v,u)$ between $u$ and $v$ such that
\Eq{CMV}{
  \frac{r(u)-r(v)+r'(v)(v-u)}{(v-u)^2}
  =\frac{r''(w)(w-u)}{2(w-u)}=\frac12r''(w).
}

Consider first the case when $K$ is not a singleton. Then there exist at least two distinct elements of $K$, say $u_1$ and $u_2$. As we have proved, $r'$ is continuously differentiable both on $J\setminus\{u_1\}$ and on $J\setminus\{u_2\}$. Therefore, it follows that $r'$ is continuously differentiable on $(J\setminus\{u_1\})\cup(J\setminus\{u_2\})=J$ and hence $r''$ is continuous on $J$. Let $u\in K$ be arbitrary. Upon taking the limit $v\to u$ in \eq{r''} and using \eq{CMV}, it follows that 
\Eq{*}{
  r''(u)=\frac12 r''(u),
}
and hence $r''(u)=0$ for all $u\in K$.

Consider now the case when $K$ is a singleton, say $K=\{u_0\}$. Then $r'$ is continuously differentiable on $J\setminus\{u_0\}$ and hence $r''$ is continuous on $J\setminus\{u_0\}$. Assume now that $r''$ has finite left and right limits at $u_0$ and denote them by $r''_-(u_0)$ and $r''_+(u_0)$, respectively. Observe that if $v<u_0$, then $v<w(v,u_0)<u_0$ and if $v>u_0$ then $v>w(v,u_0)>u_0$, therefore, upon taking the left and right limits $v\uparrow u_0$ and $v\downarrow u_0$ in \eq{r''} and using \eq{CMV} (with $u=u_0$), it follows that 
\Eq{*}{
  r_-''(u_0)=\frac12 r_-''(u_0)\qquad\mbox{and}\qquad 
  r_+''(u_0)=\frac12 r_+''(u_0).
}
These equalities imply that $r_-''(u_0)=r_+''(u_0)=0$ and hence the limit of $r''$ at $u_0$ exists and equals zero. Finally, we show that $r'$ is differentiable at $u_0$. If $v\in J\setminus\{u_0\}$, then, by the Lagrange Mean Value Theorem, there exists an element $z=z(v)$ between $v$ and $u_0$ such that 
\Eq{*}{
  \frac{r'(v)-r'(u_0)}{v-u_0}=r''(z(v)).
}
By passing the limit $v\to u_0$, we can see that $z(v)$ tends to $u_0$, and hence 
\Eq{*}{
  \lim_{v\to u_0}\frac{r'(v)-r'(u_0)}{v-u_0}
  =\lim_{v\to u_0}r''(z(v))=0.
}
This equality shows that $r'$ is differentiable at $u_0$ and $r''(u_0)=0$ and hence $r'$ is continuously differentiable at $u_0$. 
\end{proof}

One of our main results is contained in the next statement.

\Thm{Main}{Let $p:J\to\R$ and let $\varphi,\psi:J\to\R$ be continuously differentiable functions such that $\psi,p\in\CP(J)$ and $\big(\frac{\varphi}{\psi}\big)'$ is nowhere vanishing. Define $q:J\to\R $ by $q(u):=up(u)$. Then the equality \eq{2} is satisfied if and only if either $(\varphi,\psi)\sim(q,p)$, i.e., there exist four real constants $a,b,c,d$ with $ad-bc\neq0$  such that, for all $u\in J$,
\Eq{can}{
  \varphi(u)=p(u)(au+b) \qquad\mbox{and}\qquad 
  \psi(u)=p(u)(cu+d),
}
or $\varphi,\psi,p$ are infinitely many times differentiable and there exist three real constants $\alpha,\beta,\gamma$ such that \eq{cp} and
\Eq{albe0}{
  W_{q',p'}=\alpha \big(W_{q,p}\big)^3
  \qquad\mbox{and}\qquad 
  W_{\varphi',\psi'}=\beta \big(W_{\varphi,\psi}\big)^3
}
hold on $J$.}

\begin{proof} Assume that \eq{2} holds. In view of \lem{3}, we have that $p$ is continuously differentiable. By \lem{4}, it follows that the equality \eq{prod} is valid for all $(u,v)\in J^2$ and there exists a nonzero constant $\gamma\in\R$ such that \eq{cp} holds in $J$.

By \lem{Reg}, there exists an open (possibly empty) subset $H\subseteq J$ such that $\varphi,\psi,p$ are infinitely many times differentiable on $H$ and \eq{Wid} holds on $H$, where $r:=\psi/p$. If $K:=J\setminus H$ is not a singleton then $r$ is twice continuously differentiable on $J$ and $r''$ vanishes on $K$. On the other hand, if $K=\{u_0\}$ is valid for some $u_0\in J$, then $r$ is twice continuously differentiable on $J\setminus\{u_0\}$. If, additionally, $r''$ has finite left and right limits at $u_0$, then $r$ is also twice continuously differentiable at $u_0$ and $r''(u_0)=0$. 

Let $R\subseteq J$ denote the set of those points where $r$ is twice continuously differentiable. Summarizing the above possibilities, we can see that $J\setminus R$ is either empty or a singleton.

In what follows we shall distinguish two cases:

Case 1: The function $r''$ vanishes everywhere on $R$. \\ Then we show that $R=J$ holds. Indeed, if this were not true, then $J\setminus R=\{u_0\}$ for some $u_0\in J$. The function $r''$ trivially has finite left and right limits at $u_0$, hence $r$ is also twice continuously differentiable at $u_0$. This contradicts that $u_0\not\in R$ and proves that $R=J$ must be valid.
Consequently, there exist two real constants $c,d$ such that
\Eq{*}{
  r(u)=cu+d \qquad (u\in J),
}
which implies that $\psi(u)=p(u)(cu+d)$ for all $u\in J$.
Then \eq{Wd} yields that
\Eq{*}{
  \Big(\frac{\varphi}{\psi}\Big)'=\frac{\gamma p^2}{\psi^2}=\frac{\gamma}{r^2}.
}
Thus, for some real constant $\delta$ and for all $u\in J$, we get
\Eq{*}{
  \varphi(u)
  &=\psi(u)\bigg(\int\frac{\gamma}{r(u)^2}\bigg)
  =p(u)(cu+d)\bigg(\int\frac{\gamma}{(cu+d)^2}\bigg)\\
  &=p(u)(cu+d)\bigg(\delta-\frac{\gamma}{c(cu+d)}\bigg)
  =p(u)\Big(\delta(cu+d)-\frac{\gamma}{c}\Big)
  =p(u)(au+b),
}
where $a:=\delta c$ and $b:=\delta d-\gamma/c$, which also implies that $ad-bc=\gamma\neq0$. Thus \eq{can} holds for all $u\in J$.

Case II: There exists a point $u^*\in R$ such that $r''(u^*)\neq0$. \\ 
Then $u^*$ cannot belong to the set $K$, and hence $u^*\in H$. Let $U\subseteq H$ be the largest open subinterval containing $u^*$. By  \lem{Reg}, the functions $\varphi,\psi,p$ are infinitely many times differentiable on $U$, the equality \eq{Wid} is valid on $U$, and \eq{2} holds for all $(u,v)\in U^2$. Therefore, by the implication (i)$\Rightarrow$(iv) of \cite[Theorem 15]{LosPalZak21}, we have that either there exist $a,b,c,d\in\R$ with $ad-bc\neq0$ such that \eq{can} holds for all $u\in U$ or there exist two constants $\alpha,\beta\in\R$ such that
\Eq{albe}{
  \frac{W_{q',p'}}{(W_{q,p})^3}=\alpha
  \qquad\mbox{and}\qquad
  \frac{W_{\varphi',\psi'}}{(W_{\varphi,\psi})^3}=\beta
}
hold on $U$. In fact, if the first alternative were true, then $r(u)=cu+d$ implying that $r''(u)=0$ for all $u\in U$. Which cannot hold by $r''(u^*)\neq0$. Thus, the second alternative must be valid. In this case, using the identity \eq{Wid}, we have that the equality
\Eq{delta}{
\delta:=\gamma^2\beta-\alpha
=\frac{\gamma^2W_{\varphi',\psi'}}{(W_{\varphi,\psi})^3}
-\frac{W_{q',p'}}{(W_{q,p})^3}
=\frac{W_{\varphi',\psi'}-\gamma W_{q',p'}}{\gamma(W_{q,p})^3}
=-\frac{r''}{rp^4}.
}
holds on $U$. 

In what follows, we prove that $U=J$. If $U$ were a proper open subinterval of $J$, then one of the endpoints of $U$, say $w$, belongs to $J$. The point $w$ cannot belong to $H$, therefore, we must have $w\in K$. 

If $K$ is not a singleton, then $r$ is twice continuously differentiable at $w$ and $r''(w)=0$. By the continuity of $r''$ at $w$, we can see that \eq{delta} is also valid at the endpoint $w$ of $U$. Therefore, $\delta=0$, which implies that $r''/(rp^4)$ is identically zero on $U$ contradicting that $r''(u^*)\neq0$.

If $K=\{u_0\}$ for some $u_0\in J$, then $w=u_0$ and $H$ has two components:
$U$ and $V:=J\setminus(\{u_0\}\cup U)$. We may assume that $\sup U=u_0=\inf V$ (the other case is similar). From the equality \eq{delta}, we can see that $r''=-\delta rp^4$ holds on $U$, therefore, $r''$ has a finite left limit at $u_0$. On the other hand, by \lem{Reg}, the functions $\varphi,\psi,p$ are infinitely many times differentiable on $V$, the equality \eq{Wid} is valid on $V$, and \eq{2} holds for all $(u,v)\in V^2$. Therefore, by the implication (i)$\Rightarrow$(iv) of \cite[Theorem 15]{LosPalZak21}, we have that either there exist $a^*,b^*,c^*,d^*\in\R$ with $a^*d^*-b^*c^*\neq0$ such that 
\Eq{*}{
  \varphi(u)=p(u)(a^*u+b^*) \qquad\mbox{and}\qquad 
  \psi(u)=p(u)(c^*u+d^*)
}
hold for all $u\in V$ or there exist two constants $\alpha^*,\beta^*\in\R$ such that 
\Eq{*}{
  \frac{W_{q',p'}}{(W_{q,p})^3}=\alpha^*
  \qquad\mbox{and}\qquad
  \frac{W_{\varphi',\psi'}}{(W_{\varphi,\psi})^3}=\beta^*
}
hold on $V$. If the first alternative is valid, then we have that $r''=0$ is valid on $V$. If the second alternative is true, then, repeating an argument as above, we get that $r''=-\delta^* rp^4$ holds on $V$ with $\delta^*=\gamma^2\beta^*-\alpha^*$. In each case we can see that $r''$ has a finite right limit at $u_0$. Hence, in view of \lem{Reg}, the function $r$ is also twice continuously differentiable at $u_0$ and $r''(u_0)=0$. Therefore, $\delta=0$ and hence $r''$ is identically zero on $U$ contradicting that $r''(u^*)\neq0$.

The contradictions so obtained show that the endpoints of $U$ cannot belong to $J$, and hence, $U=J$. Therefore the equalities in \eq{albe} are valid on $J$ and, consequently, \eq{albe0} is also valid.

The sufficiency of the conditions of the theorem directly follows from the implication (iv)$\Rightarrow$(i) of \cite[Theorem 15]{LosPalZak21}.
\end{proof}

\Cor{EBM}{Let $p:J\to\R$ and let $\varphi,\psi:J\to\R$ be continuously differentiable functions such that $\psi,p\in\CP(J)$ and $\big(\frac{\varphi}{\psi}\big)'$ is nowhere vanishing. Define $q:J\to\R $ by $q(u):=up(u)$. Then the following assertions are equivalent:
\begin{enumerate}[(i)]
 \item The means $\B_{q,p}$ and $\B_{\varphi,\psi}$ satisfy the functional equation \eq{2}.
 \item Either $(q,p)\sim(\varphi,\psi)$ or $\varphi,\psi,p$ are twice continuously differentiable and there exist real constants $\alpha$, $\beta$ and $\gamma$ such that
 \Eq{*}{
    W_{q',p'}=\alpha\big(W_{q,p}\big)^3, \qquad
    W_{\varphi',\psi'}=\beta\big(W_{\varphi,\psi}\big)^3 \qquad\mbox{and}\qquad
    W_{\varphi,\psi}=\gamma W_{q,p}
 }
 hold on $J$.
 \item Either $(q,p)\sim(\varphi,\psi)$ or there exist real constants $a,b,c,A,B,C,\gamma$ such that
 \Eq{Q}{
   aq^2+bqp+cp^2=1, \qquad A\varphi^2+B\varphi\psi+C\psi^2=1
   \qquad\mbox{and}\qquad W_{\varphi,\psi}=\gamma W_{q,p}
 }
 hold on $J$.
 \item Either $(q,p)\sim(\varphi,\psi)$ or there exist two real polynomials $P$ and $Q$ of at most second degree which are positive on $J$ and on the range of $\varphi/\psi$, respectively, and there exist real constants $\gamma$ and $\delta$ such that
 \Eq{PQ}{
   p=\frac{1}{\sqrt{P}}, \qquad
   \psi=\frac{1}{\sqrt{Q}}\circ\frac{\varphi}{\psi}, \qquad\mbox{and}\qquad
   \bigg(\int\frac{1}{Q}\bigg)\circ\frac{\varphi}{\psi}
   =\gamma\int\frac{1}{P}+\delta.
 }
 \item Either $(q,p)\sim(\varphi,\psi)$ or there exist a continuous strictly monotone function $\omega:J\to\R$ and real constants $\alpha$ and $\beta$ such that
 \Eq{*}{
   (q,p)\sim(S_\alpha\circ\omega,C_\alpha\circ\omega)
    \qquad\mbox{and}\qquad
   (\varphi,\psi)\sim(S_\beta\circ\omega,C_\beta\circ\omega).
 }
 \item Either $(q,p)\sim(\varphi,\psi)$ or $\B_{q,p}=\B_{\omega,1}=\B_{\varphi,\psi}$ holds on $J^2$ with $\omega:=\int p^2$.
 \item Either $(q,p)\sim(\varphi,\psi)$ or there exists a continuous strictly monotone function $\omega:J\to\R$ such that  $\B_{q,p}=\B_{\omega,1}=\B_{\varphi,\psi}$ holds on $J^2$.
\end{enumerate}}

\begin{proof} The implication (i)$\Rightarrow$(ii) is a consequence of \thm{Main} because the equality \eq{cp} can be rewritten as $W_{\varphi,\psi}=\gamma W_{q,p}$.

A careful investigation of the proof of \thm{EBM} (which is the corrected form of Theorem 15 of the paper \cite{LosPalZak21}) shows that the 6th-order differentiability assumption of this theorem was only to prove the implication (i)$\Rightarrow$(ii). Therefore, the proof for the rest of the implications can be obtained following the argument of the proof of \thm{EBM} of the paper \cite{LosPalZak21}.
\end{proof}

\Cor{Last}{Let $f,g,h,k:I\to\R$ be continuously differentiable functions such that $g,k\in\CP(I)$ and the derivatives $(f/g)'$, $(h/k)'$ do not vanish in $I$. Then the assertions (i), (iii)--(vii) of \thm{EBM} are equivalent to each other. If, in addition, $f,g,h,k$ are twice continuously differentiable, then also assertion (ii) is equivalent to the other ones.}

\begin{proof} Let us define the interval $J$ and functions $p,q,\varphi,\psi:J\to\R$ as in \lem{fghk}. Then $p,q$ and $\varphi,\psi$ are continuously differentiable such that $\psi,p\in\CP(J)$ and $\big(\frac{\varphi}{\psi}\big)'$ is nowhere vanishing on $J$. In view of \cor{EBM}, we have the equivalence of the assertions (i)--(vii) of this result. Furthermore, by \lem{fghk}, the equation \eq{1} is equivalent to \eq{2}, that is assertion (i) of \thm{EBM} and \cor{EBM} are equivalent to each other. It is not difficult to check that the same equivalence property is valid for assertions (ii)--(vii) of these two results. The details of these computations are simple (but lengthy) therefore they are omitted.
\end{proof}


\providecommand{\bysame}{\leavevmode\hbox to3em{\hrulefill}\thinspace}
\providecommand{\MR}{\relax\ifhmode\unskip\space\fi MR }
\providecommand{\MRhref}[2]{%
  \href{http://www.ams.org/mathscinet-getitem?mr=#1}{#2}
}
\providecommand{\href}[2]{#2}

\end{document}